\documentclass[12pt]{article}
\usepackage{amssymb,latexsym,start2e,pstcol,pst-plot,amsthm,amsmath,psfig}

\def\Z{\mathbb{Z}}

\def\F{\mathcal{F}}
\def\<{\left<}
\def\>{\right>}
\def\ds{\displaystyle}

\begin{document}
\pagestyle{plain}

\title{What power of two divides a weighted Catalan number?
}
\author{Alexander Postnikov\\[-5pt]
\small Department of Mathematics, Massachusetts Institute of Technology,\\[-5pt]
\small Cambridge, MA 02139, \texttt{apost@math.mit.edu}\\[5pt]
Bruce E. Sagan\\[-5pt]
\small Department of Mathematics, Michigan State University,\\[-5pt]
\small East Lansing, MI 48824, \texttt{sagan@math.msu.edu}
}

\date{January 10, 2006
        \\[10pt]
	\begin{flushleft}
	                                       \\[5pt]
	\end{flushleft}}

\maketitle

\begin{abstract}
Given a sequence of integers $b=(b_0,b_1,b_2,\ldots)$ one gives a Dyck
path $P$ of length $2n$ the weight 
$$
\wt(P) = b_{h_1} b_{h_2}\cdots b_{h_n},
$$
where $h_i$ is the height of the $i$th ascent of $P$.  The corresponding
weighted Catalan number is
$$
C_n^b=\sum_P \wt(P),
$$
where the sum is over all Dyck paths of length $2n$.  So, in
particular, the ordinary Catalan numbers $C_n$ correspond to $b_i=1$ for all
$i\ge0$.  Let $\xi(n)$ stand for the base two exponent of $n$, i.e.,
the largest power of 2 dividing $n$.
We give a condition on $b$ which implies that
$\xi(C_n^b)=\xi(C_n)$. In the special case $b_i=(2i+1)^2$, this settles
a conjecture of Postnikov about the number of plane Morse links.  Our proof
generalizes the recent combinatorial proof of Deutsch
and Sagan of the classical formula for $\xi(C_n)$.
\end{abstract}

\section{Introduction}

The {\it Catalan numbers\/}
$$
C_n = \frac{1}{n+1}\binom{2n}{n}
$$ 
have many interesting arithmetic properties.  For example, the
following result which essentially dates back to Kummer (see Dickson's
book~\cite{dic:htnI} for details) describes their divisibility by powers 
of 2.  Let $s(n)$ be the sum of digits in the binary expansion of $n$.
Also, let $\xi(n)$ denote the base two {\it exponent\/} of $n$, i.e.,
the largest power of two dividing $n$.

\bth
\label{kummer}
\label{th:catalan_mod_2}
We have\\
\eqed{\xi(C_n)=s(n+1)-1.}
\eth

A combinatorial proof of this result was recently given by Deutsch and
Sagan~\cite{ds:ccm} using group actions.
In this paper, we extend this result and their proof to various
weighted Catalan numbers. 
Note that, unlike the Catalan numbers, weighted Catalan numbers need
not have simple multiplicative formulas.  So determining their divisibility
properties is more subtle than for the usual Catalan numbers.

Let $b=(b_0,b_1,b_2,\dots)$ be a fixed infinite sequence of integers.
Define the {\it weighted Catalan numbers\/}, $C_n^b$,
as the coefficients of the expansion of the continued fraction: 
\beq
\label{Cnb}
\frac{\ds 1}{\ds 1-\frac{\ds b_0\,x}{\ds 1-\frac{\ds b_1\,x}
{\ds 1- \frac{b_2\,x}{\ds 1-\frac {\ds b_3\, x}{\ds 1 - \cdots } }  }}}
= \sum_{n\geq 0} C_n^b \,x^n.
\eeq
If $b=(1,1,1,\dots)$, then $C_n^b$ is the usual Catalan number $C_n$.

Combinatorially, the $C_n^b$ count Dyck paths with certain weights.
Recall that a {\it Dyck path\/} $P$ of length $2n$ 
is a sequence of points in the upper half-plane of the integer lattice
$$
(x_0,y_0) = (0,0),\ (x_1,y_1),\ \ldots,\  (x_{2n},y_{2n}) = (2n,0),
$$
such that each step $s_i=[x_i-x_{i-1},y_i-y_{i-1}]$ has the form $[1,1]$ or $[1,-1]$.
Let us say that  step $s_i$ has {\it height\/} 
$y_{i-1}$.  Define the {\it weight\/} of a Dyck path $P$ to be
the product 
$$
\wt(P) =b_{h_1}b_{h_2}\cdots b_{h_n},
$$
where $h_1,\dots,h_n$ are the heights of its steps of the form $[1,1]$.
Then the following proposition is well known and easy to prove, i.e.,
see the book of Goulden and Jackson~\cite[Ch.~5]{gj:ce}.
\bpr
\label{dyck}
We have
$$
C_n^b=\sum_{P} \wt(P),
$$
where the sum is over all Dyck paths of length $2n$.\hqed
\epr

For example, we have 
$$
C_3^b = b_0\,b_0\,b_0 + b_0\, b_0\,b_1+b_0\,b_1\,b_0 + b_0\,b_1\,b_1 + b_0\,b_1\,b_2,
$$
where the five terms correspond to the five Dyck paths of length 6.
As another example, 
if $b=(1,q,q^2,q^3,\dots)$, then the weighted Catalan number $C_n^b$ is
equal to the {\it $q$-Catalan number\/} 
$$
C_n(q) = \sum_P q^{\mathrm{area}(P)},
$$ 
where the sum is 
over Dyck paths $P$ of length $2n$ and $\mathrm{area}(P)$ denotes
the area between $P$ and the lowest possible path. 
The continued fraction~\ree{Cnb} with $b_i = q^i$ is know as
the {\it Ramanujan continued fraction}.

\begin{figure}
\begin{tabular}{cc}
\hspace{40pt}
\psfig{figure=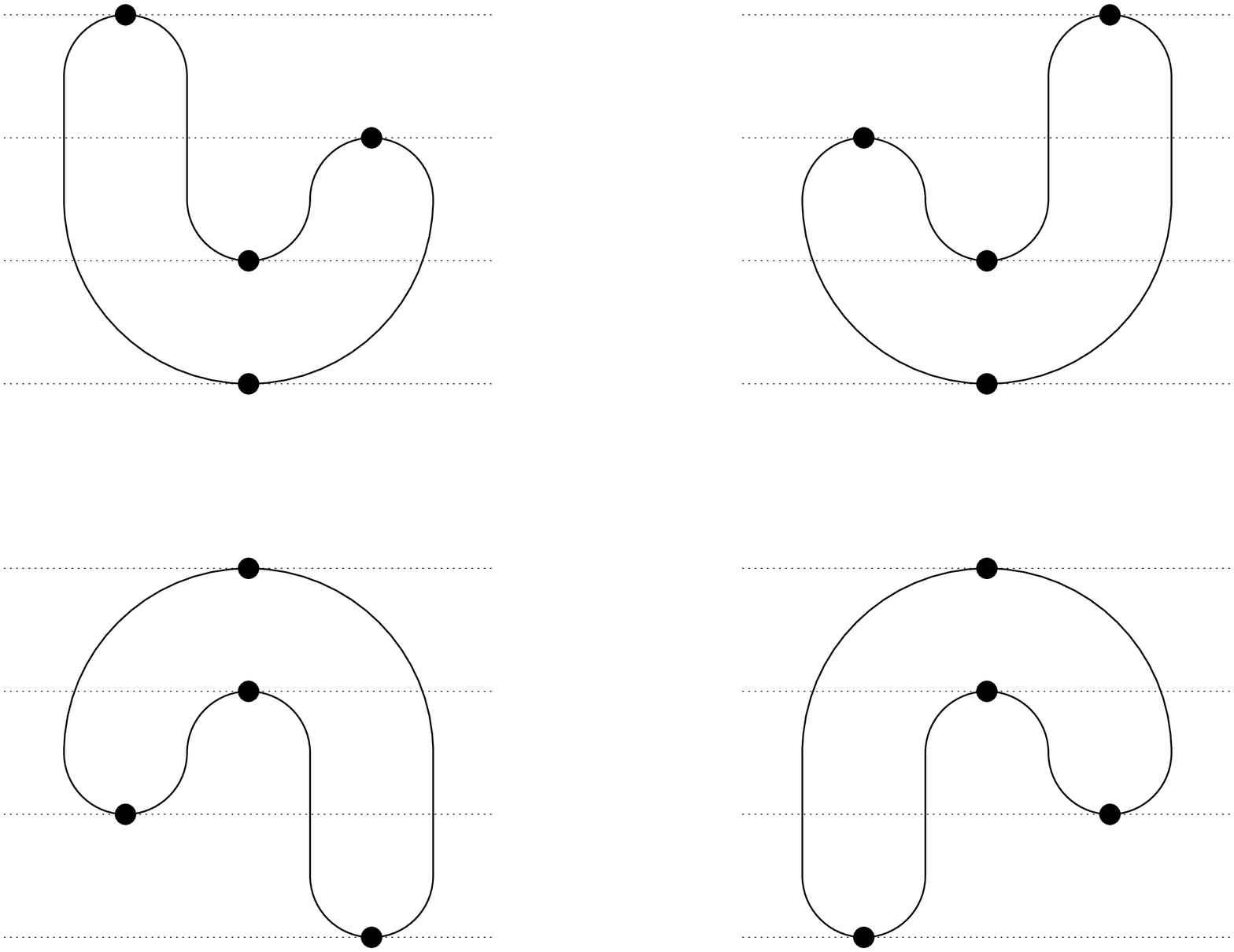,height=1in} 
&
\hspace{40pt}
\psfig{figure=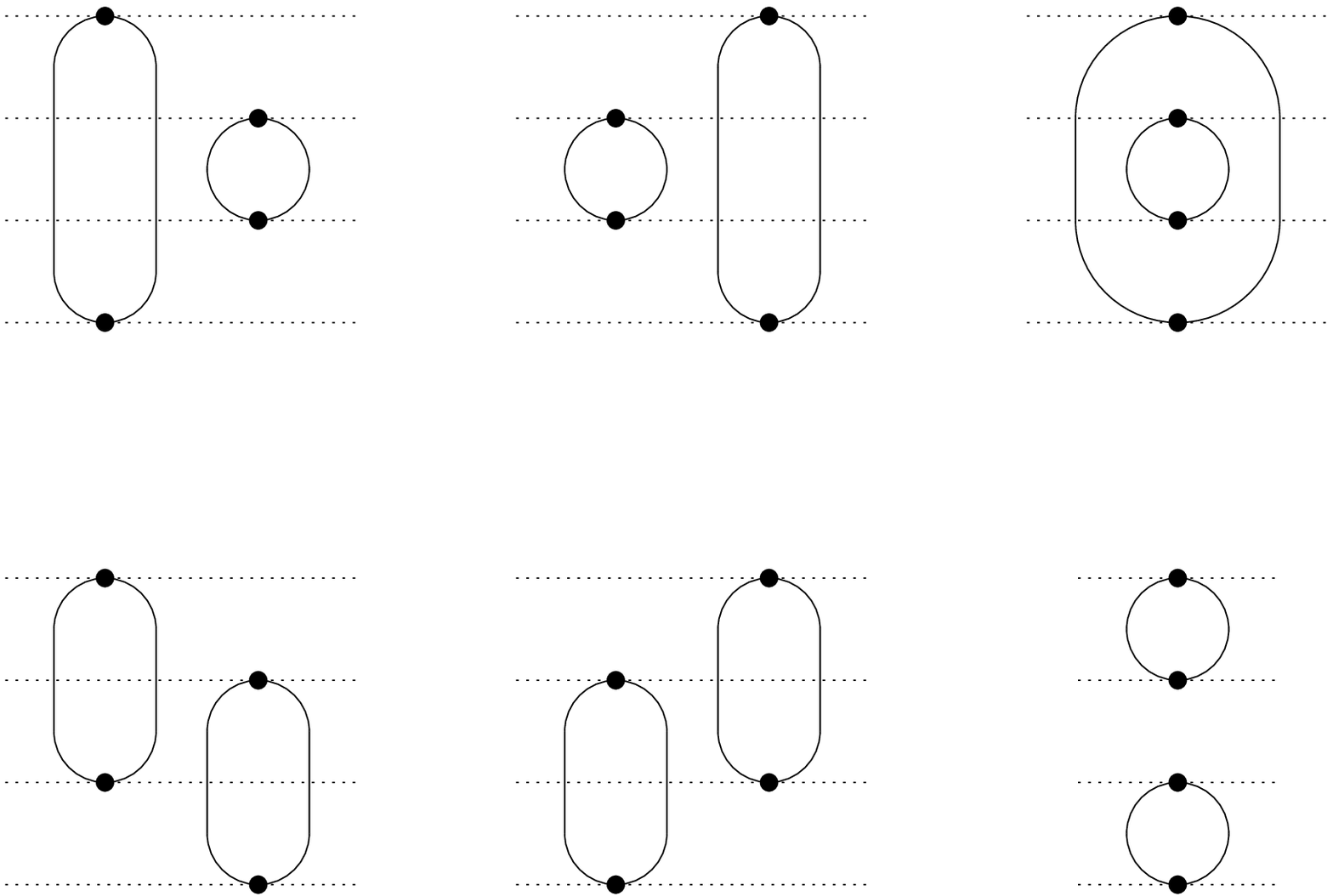,height=1in} \\[10pt]
\hspace{35pt}
(a) Connected links (curves)
&
\hspace{35pt}
(b) Disconnected links
\end{tabular}
\caption{All plane Morse links of order 2}
\label{links}
\end{figure}

Our main result, Theorem~\ref{main} below, gives a necessary condition 
on the sequence $b$ so that
$$
\xi(C_n^b)=\xi(C_n)=s(n+1)-1.
$$
As a special case, we obtain a conjecture of Postnikov~\cite{pos:cmc}
about plane Morse links.  A {\it plane Morse curve\/} is a simple
curve $f:S^1\ra \bbR^2$ (i.e., a smooth injective map) such that, for
the height function $h:(x,y)\ra y$, the map 
$h\circ f$ has a finite number of isolated nondegenerate critical
points with distinct values.  All curves are oriented clockwise.
See, for example, Figure~\ref{links}.  As one goes around a Morse curve,
the sequence formed by the critical values is alternating and returns
to where it started.  So the number of critical values must be an even
integer $2n$, and $n$ is called the {\it order\/} of the curve.  The
{\it combinatorial type\/} of a Morse curve is its connected component
in the space of all  Morse curves.  So Figure~1(a) depicts the four
plane Morse curves of order 2 up to combinatorial type.

A {\it plane Morse link\/} is a disjoint union of plane Morse curves.  All our
definitions for Morse curves carry over in the natural way to links.
In particular, the order of a link is the sum of the orders of its
components.   Figure~1(b) shows the six disconnected Morse links of
order 2, again up to combinatorial type.  Let $L_n$ be the number of
combinatorial types of plane
Morse links of order $n$.  Then we have just seen that $L_2=4+6=10$.
The connection with weighted Catalan numbers is made by the following
theorem.
\bth[\cite{pos:cmc}]  The numbers $L_n$ satisfy
$$
\sum_{n\geq 0} L_n \,x^n
=
\frac{\ds 1}{\ds 1-\frac{\ds 1^2\,x}{\ds 1-\frac{\ds 3^2\,x}
{\ds 1- \frac{5^2\,x}{\ds 1-\frac {\ds 7^2\, x}{\cdots } }  }}}
$$
and so\\
\eqed{L_n=C_n^{(1^2,\ 3^2,\ 5^2,\ 7^2,\ \ldots)}.}
\eth

An easy corollary of our main theorem will be a proof of 
Conjecture~3.1 from~\cite{pos:cmc}.  It can also be
found listed as Problem~6.C5(c) in the Catalan Addendum to the second volume of
Stanley's {\it Enumerative Combinatorics\/}~\cite{sta:ca}.
\bcon[\cite{pos:cmc,sta:ca}]  
\label{pos}
We have\\
\eqed{\xi(L_n)=\xi(C_n)=s(n+1)-1.}
\econ

\section{The main theorem}

Let $\bbZ_{\ge0}$ denote the nonnegative integers.  The {\it difference
operator\/}, $\De$, acts on functions $f:\bbZ_{\ge0}\ra\bbZ$ by
$$
(\De f)(x)=f(x+1)-f(x).
$$
Note that we can regard our sequence $b$ as such a function where
$b(x)=b_x$.  We can now state our main result, using $c\ |\ d$ as usual to
mean that $c$ divides evenly into $d$.
\bth
\label{main}
Assume that the sequence $b$ satisfies
\ben
\item $b(0)$ is odd, and
\item $2^{n+1}\ |\ (\De^n b)(x)$ for all $n\ge 1$ and $x\in\bbZ_{\ge0}$.
\een
Then
$$
\xi(C_n^b)=\xi(C_n)=s(n+1)-1.
$$
\eth

To prove this, we will also have to consider the {\it shift
operator\/}, $S$, acting on functions $f:\bbZ_{\ge0}\ra\bbZ$ by
$$
(Sf)(x)=f(x+1).
$$
It is well known and easy to verify that we have the product rule
$$
\Delta(f\cdot g) = \Delta(f)\cdot g + S(f)\cdot \Delta(g).
$$
which generalizes to
\begin{equation}
\label{eq:Delta^n}
\Delta^n(f\cdot g) = 
\sum_{k=0}^n {n\choose k} \Delta^{n-k}(S^k(f))\cdot \Delta^k(g).
\end{equation}

Let $\F$ be the set of functions $f:\bbZ_{\ge0} \to\Z$ such that 
\begin{enumerate}
\item[(a)] $f(x)$ is odd for all $x\in\bbZ_{\ge0}$, and
\item[(b)] $2^{n+1}\ |\ (\Delta^n f)(x)$ for all $n\ge 1$ and $x\in\bbZ_{\ge0}$.
\end{enumerate}
Note that because of (b), we can replace (a) by the seemingly weaker
condition that $f(0)$ is odd.  We will need the following lemma.
\ble
\label{closed}
The set $\F$ is closed under the three operations
$$
f\mapsto Sf,\quad
(f,g)\mapsto f\cdot g, \qmq{and}
(f,g)\mapsto \<f, g\> := \frac{f(x+1)\, g(x) + f(x)\, g(x+1)}{2}.
$$
\ele
\prf
Closure under $S$ is obvious.  Now
assume that $f(x),\, g(x) \in \F$.  Then the value $(f\cdot g)(x)$ 
is clearly odd.  And~(\ref{eq:Delta^n}) shows that the divisibility
criterion is satisfied.  Thus $f\cdot g\in \F$.

We can write the second operation as 
$$
\<f,g\> = f\cdot g + \frac{\Delta(f)\cdot g+f\cdot \Delta(g)}{2}.
$$
Since $\Delta(f)$ and $\Delta(g)$ are divisible by $4$ and $(f\cdot g)(x)$ is 
odd, we deduce that $\<f,g\>(x)$ is odd.
By~(\ref{eq:Delta^n}), 
$\Delta^n(\Delta(f)\cdot g)$ and  $\Delta^n(f\cdot \Delta(g))$ are
divisible by $2^{n+2}$.   Thus $\<f,g\>\in \F$.
\hqedm

Deutsch and Sagan used the interpretation of $C_n$ in terms of binary
trees to prove Theorem~\ref{kummer}.  So we will need to review this
method and translate our weight function into this setting.  A {\it
binary tree\/}, $T$, is a rooted tree where every vertex has a right
child, a left child, both children, or no children.  We also consider the
empty tree to be a binary tree.  Let $\cT_n$ be the set of binary
trees with $n$ vertices.  Then  one of the standard interpretations
of the Catalan numbers is that 
\beq
\label{Tn}
C_n=\#\cT_n.
\eeq

Let $G_n$ be the group of symmetries of the binary tree which is complete to 
depth $n$ (having all of its leaves at distance $n$ from the root).
The group $G_n$ is generated by reflections which exchange
the left and right subtrees associated with a vertex.  
Then $G_n$ acts on $\cT_n$ with two trees being in the same orbit if
they are isomorphic as rooted trees if we forget about the information
concerning left and right children.
Deutsch and Sagan show in~\cite[Section 2]{ds:ccm} that $\#G_n$ is a power of 2,
so the cardinality of any $G_n$-orbit is as well.
They also show that the minimal size of a $G_n$-orbit is
$2^s$ where $s=s(n+1)-1$.  Moreover, orbits of the minimal size can be
identified with binary total partitions on the set $\{1,2,\ldots,s\}$, whose number
is $(2s-1)!!=1\cdot 3\cdot 5\cdots(2s-1)$ as shown by
Schr\"oder~\cite{sch:vcp}, see also~\cite[Example 5.2.6]{sta:ec2}.  For ease of reference, we
summarize these facts in the following lemma.
\ble[\cite{ds:ccm}]
\label{DS}
Let $\cO$ be an orbit of $G_n$ acting on $\cT_n$ and let $s=s(n+1)-1$.
Then
\ben
\item $\#\cO=2^t$ for some $t\ge s$, and
\item $\#\cO=2^s$ for exactly $(2s-1)!!$ orbits.\hqed
\een
\ele

It is now easy to prove Theorem~\ref{kummer} using this lemma and
equation~\ree{Tn}.  To generalize the proof, consider any fixed
function $b\in\cF$ and define the corresponding {\it weight of a binary
tree $T$\/} to be the function
$$
w_b(T) = w_b(T;x)=\prod_{v\in T} b(x+l_v),
$$
where the product is over all vertices $v$ of $T$, and $l_v$ is the
number of left edges on the unique path from the root of $T$ to $v$.
If binary
tree $T$ corresponds to Dyck path $P$ under the usual depth-first
search bijection, then it is easy to see that
\beq
\label{PT}
\wt(P)= \prod_{v\in T} b(l_v) =w_b(T;0).
\eeq

We need one last lemma for the proof of Theorem~\ref{main}.  If $\cO$
is an orbit of $G_n$ acting on $\cT_n$ then we define it's 
{\it weight\/} to be
\beq
\label{wb}
w_b(\cO)=w_b(\cO;x)=\sum_{T\in\cO} w_b(T).
\eeq
\ble
For any $b\in\F$ and any orbit $\cO$ we have
$$
w_b(\cO;x) = \#\cO \cdot r_b(\cO;x),
$$
where $r_b(\cO;x)\in\cF$.
\ele
\prf
Write
\beq
\label{rb}
r_b(\cO;x)=\frac{w_b(\cO;x)}{\#\cO}.
\eeq
We induct on $n$, the number of vertices in a tree of $\cO$.  If $n=0$
then $r_b(\cO;x)=1$ for all $x$ which is clearly in $\cF$.

Now suppose $n\ge1$ so that $\cO$ contains a nonempty tree $T$.
Let $T_1$ be the subtree of $T$ consisting of the left child of the
root and all its descendents.  (So $T_1$ may be empty.)  Similarly,
define $T_2$ for the right child.
Suppose $T_1$ and $T_2$ are in orbits $\cO_1$ and $\cO_2$,
respectively.  If $\cO_1=\cO_2$ then
$\#\cO = \#\cO_1\cdot\#\cO_2$ and
$$
w_b(T) = b(x) \cdot w_{S(b)}(T_1)\cdot w_b(T_2).
$$
If $\cO_1\neq\cO_2$, then
$\#\cO = 2\cdot \#\cO_1\cdot \#\cO_2$ and
$$
w_b(T) = b(x)\ 
[w_{S(b)}(T_1)\cdot w_b(T_2)+w_b(T_1)\cdot w_{S(b)}(T_2)].
$$
In the both cases, it follows from equations~\ree{wb} and~\ree{rb} that
$$
r_b(\cO;x) = b(x)\cdot 
\frac{r_b(\cO_1;x+1)\cdot r_b(\cO_2,x)+r_b(\cO_1;x)\cdot r_b(\cO_2;x+1)}{2}.
$$
So by Lemma~\ref{closed} and induction we have $r_b(\cO;x)\in\cF$ as desired.
\hqedm

{\bf Proof of Theorem~\ref{main}.\hspace{7pt}}  Combining
Proposition~\ref{dyck}, the previous lemma, and equations~\ree{PT}
and~\ree{wb} gives
\beq
\label{sum}
C_n^b=\sum_{P}\wt(P)=\sum_{T\in\cT_n} w_b(T;0)
=\sum_{\cO} \#\cO\cdot r_b(\cO;0)
\eeq
where the integers $r_b(\cO;0)$ are all odd since $r_b(\cO)\in\cF$.
So $\xi(\#\cO\cdot r_b(\cO;0))=\xi(\#\cO)$ for all orbits $\cO$.  It
now follows from Lemma~\ref{DS} that $\xi=s$ for an odd number (namely
$(2s-1)!!$) of summands in the last summation in~\ree{sum} and that $\xi>s$
for the rest.  We conclude that $\xi(C_n^b)=s=\xi(C_n)$.\hqedm

As corollaries, we can prove Conjecture~\ref{pos} and give
information about divisibility of the $q$-Catalan numbers.
\bco
1.
The number of combinatorial types of plane Morse links of order $n$
satisfies
$$
\xi(L_n)=\xi(C_n)=s(n+1)-1.
$$
2.
If $q\Cong 1\ (\Mod 4)$ then the $q$-Catalan numbers satisfy
$$
\xi(C_n(q))=\xi(C_n)=s(n+1)-1.
$$
\eco
\prf
For the first assertion, it suffices to show that the function
$b(x)=(2x+1)^2$ is in $\cF$.  Clearly $b(x)$ is odd for all
$x\in\bbZ_{\ge0}$.  Furthermore $\De b=8(x+1)$, $\De^2 b=8$, and $\De^n b=0$
for $n\ge3$.  So the divisibility condition also holds.

For the second statement, we need the function $b(x)=q^x$ to be in
$\cF$.  Since $q$ is odd, so is $b(x)$.  Also $\De^n b=(q-1)^n q^x$
for $n\ge1$.  So the hypothesis on $q$ implies that $\De^n b$ is
divisible by $4^n=2^{2n}$ which is more than needed.
\hqed

\section{An open problem}

Consider the {\it Catalan sequence\/}
$$
C=(C_0,C_1,C_2,\ldots).
$$
Theorem~\ref{kummer} implies immediately that $C_n$ is odd if and
only if $n=2^k-1$ for some $k\ge0$.  It follows that
the $k$th block of
zeros in the sequence $C$ taken modulo 2 has length $2^k-1$
(where we start numbering with the first block).  Alter and
Kubota~\cite{ak:ppp} have generalized this result to arbitrary primes
and prime powers.  One of their main theorems is as follows.
\bth[\cite{ak:ppp}]
Let $p\ge3$ be a prime and let $q=(p+1)/2$. The length of the $k$th
block of zeros in $C$ modulo $p$ is
$$
\frac{p^{\xi_q(k)+\de_{3,p}+1}-3}{2}
$$
where $\xi_q(k)$ is the largest power of $q$ dividing $k$ and
$\de_{3,p}$ is the Kronecker delta.\hqedm 
\eth

Deutsch and Sagan~\cite{ds:ccm} have improved on this theorem when
$p=3$ by giving a complete characterization of the residues in $C$
modulo 3.  However, the demonstrations of all these results rely
heavily on the expression for $C_n$ as a product.  It would be
interesting to find analogous theorems for $C_n^b$, but new proof
techniques would have to be found.

\bigskip
\bibliographystyle{plain}
\begin{small}
\bibliography{ref}
\end{small}

\end{document}